\documentclass[smallextended]{svjour3}
\usepackage{amsfonts}
\usepackage{calrsfs}
\usepackage{mathptmx}

\newcommand\Prob{{\mathbb P}}
\newcommand\E{{\mathbb E}}
\newcommand\R{{\mathbb R}}

\newcommand\N{{\mathbb N}}
\newcommand\Z{{\mathbb Z}}
\newcommand\I{{\mathbb I}}
\newcommand\Var{{\mathbb V{\rm ar}}}

\smartqed
\newtheorem{Theorem}{Theorem}
\newtheorem{Lemma}[Theorem]{Lemma}
\newtheorem{Corollary}[Theorem]{Corollary}
\newtheorem{Proposition}[Theorem]{Proposition}

\begin{document}
\title{Strong law of large numbers for a function\\ of
the local times of a transient random walk in $\Z^d$}
\titlerunning{Law of Large Numbers for Local Times}
\author{I. M. Asymont \and D. Korshunov (*)}

\institute{I. M. Asymont \at
              Financial University under the Government 
of the Russian Federation, Russia\\
              \email{i.asymont@mail.ru}           
           \and
           D. Korshunov (*) \at
              Lancaster University, UK\\
\email{d.korshunov@lancaster.ac.uk}
}

\maketitle

\begin{abstract}
For an arbitrary transient random walk $(S_n)_{n\ge 0}$ in $\Z^d$, 
$d\ge 1$, we prove a strong law of large numbers for the spatial sum
$\sum_{x\in\Z^d}f(l(n,x))$
of a function $f$ of the local times $l(n,x)=\sum_{i=0}^n\I\{S_i=x\}$. 
Particular cases are 
the number of

(a) visited sites (first time considered by 
Dvoretzky and Erd\H{o}s in \cite{DE}),
which corresponds to a function $f(i)=\I\{i\ge 1\}$; 

(b) $\alpha$-fold self-intersections of the random walk
(studied by Becker and K\"{o}nig in \cite{BK}),
which corresponds to $f(i)=i^\alpha$; 

(c) sites visited by the random walk exactly $j$ times
(considered by Erd\H{o}s and Taylor in \cite{ET} 
and by Pitt \cite{Pitt}), where $f(i)=\I\{i=j\}$.


\end{abstract}

\section{Introduction and main results}

Let $X_1$, $X_2$, \ldots be a sequence of independent identically
distributed random vectors valued in $\Z^d$, $d\ge 1$.
Consider a {\it random walk} generated by $X_n$'s, $S_0:=0$, 
$S_n : =X_1+\ldots+X_n$, and the number of visits to a 
site $x\in{\mathbb{Z}^d}$ up to time $n$
which is called the {\it local time} of $x$,
\begin{equation}\label{f1.2}
l(n,x) : =\sum_{i=0}^n \I\{S_i=x\}.
\end{equation}
Define random variables
\begin{equation}\label{f1.3}
L_n(\alpha)\ :=\ \sum_{x\in\Z^d:\ l(n,x)>0} l^\alpha(n,x), \quad \alpha\ge0.
\end{equation}
In particular, the $L_n(0)=|\{S_0,...,S_n\}|$  represents the number 
of distinct sites visited by the random walk up to time $n$,
called the {\it range} of $(S_n)_{n\ge 1}$. 
The case $\alpha=1$ is trivial because $L_n(1)=n+1$. 
The value of $L_n(2)$ is the number of so called {\it self-intersections}
of a random walk. For an integer $\alpha$ the value of $L_n(\alpha)$ 
is the number of $\alpha$-fold self-intersections up to time $n$.

It is known that for a {\it recurrent} random walk 
the quotient $L_n(0)/n$ tends to $0$ as $n\to\infty$ 
(see, e.g. Spitzer \cite[Ch. 1, Sect. 4, Theorem 1]{S}),
which assumes a slower growing normalising sequence for a proper limit
in the law of large numbers.
As shown in Dvoretzky and Erd\H{o}s \cite[Theorem 3]{DE} for a simple random walk
and in \v{C}ern\'y \cite{Cerny} for a general one with zero drift
and finite covariance matrix, it is $n/\log n$ in 2 dimensions. 

In present article we show that the law of large numbers for $L_n(\alpha)$
with a non-zero limit and normalising sequence $n$ holds
in any dimension $d$ for any {\it transient} random walk, 
that is, when the probability of its return to the origin, 
$$
\gamma\ :=\ \Prob\{S_n\neq 0\mbox{ for all }n\ge 1\},
$$ 
is strictly positive, $\gamma>0$. We assume in addition that 
$\gamma<1$ which excludes a trivial case where either $l(n,x)=1$ or $0$ 
for all $x$ with probability $1$, and hence $L_n(\alpha)=n+1$.

The result for $L_n(\alpha)$ we are interested in follows from 
the following more general result.
Consider a function $f:\Z^+\to\R$ and a spatial sum
$$
G_n(f)\ :=\ \sum_{x\in\Z^d}f(l(n,x)).
$$
In particular, for a power function $f(i)=i^\alpha$,
we get $L_n(\alpha)=G_n(f)$.

\begin{Theorem}\label{t1.2}
Let the random walk $(S_n)_{n\ge 0}$ be transient and $f:\Z^+\to\R$ 
be a function satisfying 
\begin{eqnarray}\label{cond.on.f}
\sum_{j=1}^\infty f^2(j)j(1-\gamma)^j &<& \infty.
\end{eqnarray}
Then
\begin{eqnarray}\label{Gf.conv.L21}
\frac{G_n(f)}{n} &\to&
\gamma^2\sum_{j=1}^\infty f(j)(1-\gamma)^{j-1}
\quad \mbox{as }n\to\infty
\end{eqnarray}
in mean square and with probability $1$.
\end{Theorem}

The proof of Theorem \ref{t1.2} is given in Section \ref{sec:proof1}. 
In Sections \ref{sec:expectation} and \ref{sec:variance} 
we discuss an asymptotic behaviour of the expectation 
and variance of $G_n(f)$ as $n\to\infty$ respectively,
needed further in the proofs.

The following corollaries are immediate.

\begin{Corollary}\label{t1.1}
For any $\alpha\ge 0$, it holds that
\begin{eqnarray}\label{BKresult}
\frac{L_n(\alpha)}{n} &\to&
\gamma^2\sum_{j=1}^\infty j^\alpha (1-\gamma)^{j-1}
\quad \mbox{as }n\to\infty
\end{eqnarray}
in mean square and with probability $1$.
\end{Corollary}

The case $\alpha=0$ was considered by Spitzer in \cite[Theorem 1.4.1]{S}
where convergence in probability is proven. Before then a strong
law of large numbers for $\alpha=0$ was proven for a simple random walk
by Dvoretzky and Erd\H{o}s in \cite{DE}.
In Becker and K\"{o}nig \cite{BK} the strong convergence (\ref{BKresult})
is proven for all $\alpha\ge 0$ 
(up to a gap in the proof of Proposition 2.1, 
see a comment on it in the proof of Lemma \ref{p2.1}
following equation (\ref{Q.n.j})) without any further 
conditions in the case $d\ge 3$, however in the cases
$d\in{\{1,2\}}$ it is assumed there that either the steps 
$X_i$ are square integrable or, for some $\eta>0$ and $C<\infty$,
\begin{eqnarray}\label{bound.for.P.Sn=0.pre}
\sum_{k=n}^\infty \Prob\{S_k=0\} &\le& \frac{C}{n^\eta}\quad\mbox{for all }n. 
\end{eqnarray}

\begin{Corollary}\label{cor:j}
Let 
$J\subset\N$. Then, with probability $1$,
\begin{eqnarray*}
\frac{1}{n} \sum_{x\in\Z^d} \I\{l(n,x)\in J\} &\to&
\sum_{j\in J}\gamma^2(1-\gamma)^{j-1} \quad \mbox{as }n\to\infty.
\end{eqnarray*}
\end{Corollary}

If $J$ is a singleton $\{j\}$, then we get the strong law of large numbers
for the number of sites visited exactly $j$ times up to time $n$.
For these statistics, the last corollary generalises Theorem 12 
in Erd\H{o}s and Taylor \cite{ET} from a simple random walk in 
$d\ge 3$ dimensions to an arbitrary transient random walk;
a general result for transient random walks on a countable Abelian group
was proven by induction on $j$ by Pitt in \cite{Pitt}.
Notice that, for an arbitrary $J$, say $J$ the set of all odd numbers,
Corollary \ref{cor:j} can be reduced to the singleton case,
once we know the strong law of large numbers for the range of $S_n$.

The growth condition (\ref{cond.on.f}) is satisfied for all
subexponential functions $f(i)$ of order $e^{o(i)}$ as $i\to\infty$, 
and also for exponentially growing functions of order $O(e^{ci})$
with exponent coefficient $c<\lambda_*/2$ where
$$
\lambda_*:=\log\frac{1}{1-\gamma}.
$$
It is very likely that the condition (\ref{cond.on.f}) may be relaxed to the 
condition (\ref{cond.on.f.w}) below because under the latter condition we have 
\begin{eqnarray*}
\E |f(l(\infty,0))| &=& \sum_{k=1}^\infty 
|f(k)|(1-\gamma)^{k-1}\gamma\ <\ \infty,
\end{eqnarray*}
and since the number of visited sites up to time $n$ is not greater than $n$,
it clearly indicates that the family $\{G_n(f)/n,n\ge 1\}$ 
is stochastically bounded. But if we only assume (\ref{cond.on.f.w}),
then it requires a much more delicate analysis
compared to the estimation of the variance carried out in Lemma \ref{l:var},
as it happens when we prove a strong law of large numbers
for a random walk where existence of the second moment 
of jumps essentially simplifies proving technique.
In the result below we show how it can be done under some additional
technical assumptions. 

\begin{Theorem}\label{thm:exp.conv}
Let, for some $C<\infty$ and $\varepsilon>0$,

{\rm(i)} either the condition 
\begin{eqnarray}\label{bound.for.P.Sn=0}
\sum_{k=1}^n k\Prob\{S_k=0\} &\le& Cn^{1-\eta}\quad\mbox{for all }n
\end{eqnarray}
hold for some $\eta\in(0,1)$ and 
$|f(i)|\le Ce^{i\lambda_*}/i^{2+\varepsilon}$ for all $i$,

{\rm(ii)} or the condition 
\begin{eqnarray}\label{bound.for.P.Sn.conv}
\sum_{k=1}^n k\Prob\{S_k=0\} &\le& C\log n\quad\mbox{for all }n
\end{eqnarray}
hold and $|f(i)|\le Ce^{i\lambda_*}/i\log^{2+\varepsilon}i$ for all $i>0$.

Then the convergence (\ref{Gf.conv.L21}) holds with probability $1$.
\end{Theorem}

For the proof, see Section \ref{sec:max}. 
It is based on truncation technique and on a strong limit theorem 
for the maximal local time, $l(n):=\max\{l(n,x),\ x\in\Z^d\}$, 
see Proposition \ref{l:ln.limit.thm} there. 

Notice that the condition (\ref{bound.for.P.Sn=0}) is equivalent to
(\ref{bound.for.P.Sn=0.pre}). Indeed, on the one hand,
it follows from (\ref{bound.for.P.Sn=0.pre}) that
\begin{eqnarray*}
\sum_{k=1}^n k\Prob\{S_k=0\} &=& \sum_{j=1}^n \sum_{k=j}^n\Prob\{S_k=0\}
\ \le\ C\sum_{j=1}^n j^{-\eta}\ \le\ \frac{C}{1-\eta}n^{1-\eta}.
\end{eqnarray*}
On the other hand, it follows from (\ref{bound.for.P.Sn=0}) that, for all $m$,
\begin{eqnarray*}
C(2m)^{1-\eta}\ \ge\ \sum_{k=m}^{2m-1} k\Prob\{S_k=0\} &\ge& 
m\sum_{k=m}^{2m-1} \Prob\{S_k=0\},
\end{eqnarray*}
hence
\begin{eqnarray*}
\sum_{k=n}^\infty \Prob\{S_k=0\} &=& 
\sum_{j=0}^\infty \sum_{k=n2^j}^{n2^{j+1}-1} \Prob\{S_k=0\}
\ \le\ C2^{1-\eta}\sum_{j=0}^\infty (n2^j)^{-\eta}
\ =\ \frac{2C}{2^\eta-1}n^{-\eta}.
\end{eqnarray*}
Also notice that, in $d\ge 3$ dimensions, if a random walk is not 
concentrated in some $3$-dimensional subspace, then the condition
(\ref{bound.for.P.Sn=0}) is valid because $\Prob\{S_n=0\}=O(1/n^{d/2})$,
due to an upper bound for the concentration function of a sum of random vectors, 
see e.g. Corollary of Theorem 6.2 in Esseen \cite{Esseen}.
For the same reason, in $d\ge 4$ dimensions, the condition
(\ref{bound.for.P.Sn.conv}) is valid for any random walk  
not concentrated in some $3$-dimensional subspace.

If the function $f$ grows faster than assumed in Theorem \ref{thm:exp.conv},
say if the condition (\ref{cond.on.f.w}) fails,
then $G_n(f)$ would require stronger normalisation than just $n$,
in order to have a proper limit as $n\to\infty$. 
The answer may be conjectured as follows:
let $\tau=\inf\{n\ge 1: S_n=S_0\}$ be the first return time to the origin, then
\begin{eqnarray*}
\E |f(l(n,0))| &=& \sum_{k=1}^n 
|f(k)|\Prob\{\widetilde\tau_1+\ldots+\widetilde\tau_{k-1}\le n\}
(1-\gamma)^{k-1}\gamma,
\end{eqnarray*}
where $\widetilde\tau_1$, $\widetilde\tau_2$, \ldots\ are independent copies of $\tau$
conditioned on $\{\tau<\infty\}$. For example, consider $f$ such that 
$f(k)\sim c_1/(1-\gamma)^k$, then
\begin{eqnarray*}
\E f(l(n,0)) &\sim& c_2\sum_{k=1}^\infty
\Prob\{\widetilde\tau_1+\ldots+\widetilde\tau_{k-1}\le n\}\quad\mbox{as }n\to\infty.
\end{eqnarray*}
As shown in \cite[Theorem 4]{DK},
in the case where $\E X_1=0$, $\E\|X_1\|^2<\infty$ and $d\ge 3$, we have 
an asymptotic relation $\Prob\{\widetilde\tau=n\}\sim c_3/n^{d/2}$ as $n\to\infty$. 

Hence, in the case $d\ge 5$, $\E\widetilde\tau_1<\infty$ and it follows 
from the renewal theorem that then $\E f(l(n,0))\sim c_2n/\E\widetilde\tau_1$,
which together with asymptotic size of the range---which is of order 
$O(n)$---indicates that the right normalisation for $G_n(f)$ should be $n^2$.

In the cases $d=3$ and $d=4$, $\E\widetilde\tau_1=\infty$ and it follows 
from Erickson's renewal theorem \cite[Theorem 5]{Erickson} 
that then $\E f(l(n,0))\sim c_4n^{1/2}$ and $c_5n/\log n$ respectively,
which in turn indicates that the right normalisation for $G_n(f)$ should be 
$n^{3/2}$ and $n^2/\log n$ respectively.

\section{Asymptotics for expectation of $G_n(f)$}
\label{sec:expectation}

In this section we discuss the asymptotic behaviour of $\E G_n(f)$ 
as $n\to\infty$. We prove the following result.

\begin{Lemma}\label{p2.1}
Let $f:\Z^+\to\R$ be a function satisfying 
\begin{eqnarray}\label{cond.on.f.w}
\sum_{j=1}^\infty |f(j)|(1-\gamma)^j &<& \infty.
\end{eqnarray}
Then
\begin{eqnarray}\label{Gf.conv}
\frac{\E G_n(f)}{n} &\to&
\gamma^2\sum_{j=1}^\infty f(j)(1-\gamma)^{j-1}
\quad\mbox{as }n\to\infty.
\end{eqnarray}
\end{Lemma}

\proof
Following Dvoretzky and Erd\H{o}s \cite{DE}, we introduce 
\begin{eqnarray*}
\gamma_n &:=& \Prob\{S_n\ne S_k\mbox{ for all } 0\le k\le n-1\},
\end{eqnarray*}
the probability that the site visited by random walk in the $n$th step
has not been visited before then; $\gamma_0=1$. As noticed in \cite{DE},
\begin{eqnarray*}
1-\gamma_n &=& \Prob\{S_n=S_k\mbox{ for some } 0\le k\le n-1\}\\
&=& \Prob\{S_{n-k}=S_0\mbox{ for some } 0\le k\le n-1\},
\end{eqnarray*}
so $\gamma_n$ equals the probability that the random walk does not return
to the origin in $n$ steps: 
\begin{eqnarray*}
\gamma_n &=& \Prob\{S_k\ne S_0\mbox{ for all } 1\le k\le n\}
\ =\ \Prob\{\tau\ge n+1\}.
\end{eqnarray*}
We observe the following monotone convergence
\begin{eqnarray}\label{f2.3}
\gamma_n-\gamma\ =\ \Prob\{n+1\le\tau<\infty\} &\downarrow& 0
\quad \mbox{as }n\to\infty.
\end{eqnarray}

Consider the following spatial sum
$$
Q_n(j)\ :=\ \sum_{x\in\Z^d} \I\{l(n,x)=j\},
$$ 
which represents the number of sites visited exactly $j$ times up to time $n$, 
hence
\begin{eqnarray}\label{decomposition}
G_n(f) &=& \sum_{j=1}^n f(j) Q_n(j).
\end{eqnarray} 
As Becker and K\"{o}nig \cite[Eq. (2.2)]{BK} do,
we use the following equality, for $j\ge 1$:
\begin{eqnarray}\label{Q.n.j}
\lefteqn{\E Q_n(j)}\nonumber\\
&=& \sum_{x\in\Z^d}{\Prob\{l(n,x)=j\}}\nonumber\\
&=& \sum_{{x\in\Z^d}\atop {0\le{k_1}<\ldots<k_j\le{n}}}
\Prob\{S_{k_1}=...=S_{k_j}=x, S_k\ne x
\mbox{ for all } k\le n,k\not\in\{k_1,\ldots,k_j\}\}\nonumber\\
&=& \sum_{0\le k_1<\ldots<k_j\le n}\gamma_{k_1}
\biggl[\prod_{i=1}^{j-1} \Prob\{\tau=k_{i+1}-k_i\}\biggr]\gamma_{n-k_j},
\end{eqnarray}
due to the Markov property of the random walk.
In \cite{BK}, the asymptotic behaviour of $\E Q_n(j)$ as $n\to\infty$ 
is argued by considering the generating function of 
$\{\E Q_n(j),n\ge 1\}$ and then referring to the Tauberian theorem,
\cite[Theorem XIII.5]{Feller}. Notice that this approach 
requires the sequence $\{\E Q_n(j),n\ge 1\}$ to be ultimately increasing 
(see Sections 1.7.3 and 1.7.4 in \cite{BGT})
which is not granted from the beginning and probably fails; 
at least such a discussion is missing in \cite{BK}. 
Notice that this problem can be fixed by first looking at the sum 
of $Q_j(n)$ over $j\ge\widetilde j$ for some $\widetilde j$ 
(this is now monotonic in $n$) and then looking at the differences.
See also Pitt \cite{Pitt} for an alternative proof.

Below we suggest another argument which does not require 
the Tauberian theorem and is only based on the transience of the random walk. 
It follows from (\ref{Q.n.j}) that
\begin{eqnarray*}
\E Q_n(j) &=& \sum_{{{n_0,n_j\ge 0}\atop 
 {n_1,\ldots,n_{j-1}\ge 1}}\atop {n_0+\ldots+n_j=n}}
 \gamma_{n_0}\gamma_{n_j}\prod_{i=1}^{j-1}\Prob\{\tau=n_i\}\\
&=& \sum_{n_0,n_j:\ 0\le n_0+n_j\le n-j+1}
 \gamma_{n_0}\gamma_{n_j}\Prob\{\tau_1+\ldots+\tau_{j-1}=n-(n_0+n_j)\},
\end{eqnarray*}
where $\tau_1$, $\tau_2$, \ldots\ are independent copies of $\tau$,
the first return time to the origin. Thus
\begin{eqnarray*}
\E Q_n(j) &=& 
\sum_{s=0}^{n-j+1}\Prob\{\tau_1+ ...+\tau_{j-1}=n-s\}
\sum_{n_0=0}^s \gamma_{n_0}\gamma_{s-n_0},
\end{eqnarray*}
hence
\begin{eqnarray}\label{Qj.represent}
\frac {\E Q_n(j)}{n} &=& \sum_{s=j-1}^n \Prob\{\tau_1+ ...+\tau_{j-1}=s\}
\frac{1}{n}\sum_{n_0=0}^{n-s} \gamma_{n_0}\gamma_{n-s-n_0}.
\end{eqnarray}
In view of the convergence (\ref{f2.3}), for any fixed $s\ge j-1$, 
\begin{eqnarray*}
\frac{1}{n}\sum_{n_0=0}^{n-s} \gamma_{n_0}\gamma_{n-s-n_0} &\to& \gamma^2
\quad\mbox{as }n\to\infty,
\end{eqnarray*}
and, moreover,
\begin{eqnarray*}
\frac{1}{n}\sum_{n_0=0}^{n-s} \gamma_{n_0}\gamma_{n-s-n_0} 
&\le& \frac{n-s+1}{n}\ \le\ 1\quad\mbox{for all }n\mbox{ and } s\ge 1.
\end{eqnarray*}
Therefore, by the dominated convergence theorem, as $n\to\infty$,
\begin{eqnarray*}
\frac {\E Q_n(j)}{n} &\to& 
\gamma^2\Prob\{\tau_1+ ...+\tau_{j-1}<\infty\}\nonumber\\
&=& \gamma^2\prod_{k=1}^{j-1}\Prob\{\tau_k<\infty\}
\ =\ \gamma^2(1-\gamma)^{j-1},
\end{eqnarray*}
owing to independence of $\tau_k$'s. In addition,
\begin{eqnarray}\label{upper.for.Qi}
\frac{\E Q_n(j)}{n} &\le& 
\Prob\{\tau_1+ ...+\tau_{j-1}<\infty\}\ =\ 
(1-\gamma)^{j-1}\quad\mbox{for all }n.
\end{eqnarray}
Then the condition (\ref{cond.on.f.w}) makes it possible to apply
dominated convergence again and to conclude that 
$$
\frac{\E G_n(f)}{n}\ =\ \sum_{j=1}^\infty f(j)\frac{\E Q_n(j)}{n}
\ \to\ \gamma^2\sum_{j=1}^\infty f(j)(1-\gamma)^{j-1}
\quad\mbox{as }n\to\infty,
$$
which completes the proof of (\ref{Gf.conv}). 
Also notice that (\ref{upper.for.Qi}) implies an upper bound
\begin{eqnarray}\label{upper.for.Q}
\E G_n(f) &\le& n\sum_{j=1}^n f(j)(1-\gamma)^{j-1}.
\end{eqnarray} 
\qed

\section{Estimation of variance of $G_n(f)$}
\label{sec:variance}

The proof of the strong law of large numbers for $L_n(\alpha)$ 
for a transient random walk given by Becker and K\"{o}nig in \cite{BK}
is based on the following upper bound for the variance of $L_n(\alpha)$: 
\begin{eqnarray*}
\Var L_n(\alpha) &\le& Cn\sum_{x\in\Z^d} \sum_{i,j=0}^n \Prob\{S_i=x\}
\Prob\{S_j=-x\}\quad\mbox{for all } n,
\end{eqnarray*}
where $C=C(\alpha)$ is a constant. 
Notice that the proof of this bound provided in \cite{BK} 
starts with an analysis of some representation for the variance of $L_n(\alpha)$,
which is only available for integer $\alpha$'s,
implication of which is necessary for further arguments 
for the strong law of large numbers for $L_n(\alpha)$
in the case of a non-integer $\alpha$.

For this reason we suggest below a different bound which works 
not only for $L_n(\alpha)$ with a non-integer $\alpha$, 
but also for $G_n(f)$ with a function $f$ other than power.
This bound provides a straightforward way for proving the strong
law of large numbers for $G_n(f)$ with $f$ satisfying the growth
condition (\ref{cond.on.f}).

\begin{Lemma}\label{l:var}
For any non-decreasing function $f$ with $f(0)=0$,
\begin{eqnarray*}
\Var G_n(f) &\le& \E G_n(f^2)
+4\sum_{i=1}^n f(i)\Delta f(i)
(1-\gamma)^{i-1} \sum_{r=1}^n r(n-r)\Prob\{S_r=0\}
\end{eqnarray*}
for all $n$ where $\Delta f(i):=f(i)-f(i-1)\ge 0$.
\end{Lemma}

\proof
In view of the representation (\ref{decomposition}),
\begin{eqnarray*}
G_n(f) &=& \sum_{x\in\Z^d}\sum_{i=1}^n f(i) \I\{l(n,x)=i\},
\end{eqnarray*}
hence
\begin{eqnarray*}
\lefteqn{\Var G_n(f)}\\
&=& 
\sum_{x,y\in\Z^d}\sum_{i,j=1}^n f(i)f(j)\Bigl(\Prob\{l(n,x)=i,l(n,y)=j\}
-\Prob\{l(n,x)=i\}\Prob\{l(n,y)=j\}\Bigr)\\
&=& \sum_{x\in\Z^d}\sum_{i=1}^n f^2(i)\Prob\{l(n,x)=i\}
+\sum_{x\not=y}\sum_{i,j=1}^n f(i)f(j)\Prob\{l(n,x)=i,l(n,y)=j\}\\
&&-\sum_{x,y\in\Z^d}\sum_{i,j=1}^n f(i)f(j)\Prob\{l(n,x)=i\}\Prob\{l(n,y)=j\},
\end{eqnarray*}
because $\Prob\{l(n,x)=i,l(n,y)=j\}=0$ if $x=y$ and $i\not= j$.
Thus, due to $f\ge 0$,
\begin{eqnarray}\label{upper.1}
\Var G_n(f) &\le& \E G_n(f^2)
+\sum_{i,j=1}^n f(i)f(j)
\Bigl(\sum_{x\not=y}\Prob\{l(n,x)=i,l(n,y)=j\}\nonumber\\
&&\hspace{35mm}-\sum_{x,y}\Prob\{l(n,x)=i\}\Prob\{l(n,y-x)=j\}\Bigr)\nonumber\\
&=:& \E G_n(f^2)+\Sigma^1_n-\Sigma^2_n,
\end{eqnarray}
and it only remains to estimate the difference of sums 
$\Sigma^1_n-\Sigma^2_n$ on the right hand side. Since $f(0)=0$, 
\begin{eqnarray*}
\lefteqn{\sum_{i,j=1}^n f(i)f(j)\Prob\{l(n,x)=i,l(n,y)=j\}}\\
&=& \sum_{i,j=1}^n \Prob\{l(n,x)=i,l(n,y)=j\}
\sum_{i_1=1}^i\sum_{j_1=1}^j \Delta f(i_1)\Delta f(j_1)\\
&=& \sum_{i_1,j_1=1}^n \Delta f(i_1)\Delta f(j_1)
\sum_{i=i_1}^n\sum_{j=j_1}^n \Prob\{l(n,x)=i,l(n,y)=j\}\\
&=& \sum_{i,j=1}^n \Delta f(i)\Delta f(j)\Prob\{l(n,x)\ge i,l(n,y)\ge j\},
\end{eqnarray*}
and similar equalities hold for ordinary sums. Therefore,
\begin{eqnarray*}
\Sigma^1_n-\Sigma^2_n &=& \sum_{i,j=1}^n \Delta f(i)\Delta f(j)
\biggl(\sum_{x\not=y}\Prob\{l(n,x)\ge i,l(n,y)\ge j\}\\
&&\hspace{40mm}
-\sum_{x,y}\Prob\{l(n,x)\ge i\}\Prob\{l(n,y)\ge j\}\biggr),
\end{eqnarray*}
where $\Delta f(i)\ge 0$ for all $i$ because $f$ is non-decreasing,
and the tail probabilities do not decrease as $n$ grows, 
which makes it possible to perform a required analysis of the double sum.
Let us decompose the event $B=B(x,y,i,j):=\{l(n,x)\ge i,l(n,y)\ge j\}$ 
for $x\not= y$ as a union of four disjoint events 
$B\cap B_{xy}$, $B\cap B_{yx}$, $B\cap B_{xyx}$ and $B\cap B_{yxy}$, where
\begin{eqnarray*}
B_{xyx} &:=& \{S_{n_1}=x, S_{n_2}=y,S_{n_3}=x
\mbox{ for some }n_1<n_2<n_3\le n\},\\
B_{xy} &:=& \{\mbox{all visits to }x\mbox{ occur before all visits to }y\}.
\end{eqnarray*}
Denote by $\tau_x(i)$ the time of $i$th visit to $x$ 
by the random walk $(S_n)_{n\ge 0}$. Then the event $B\cap B_{xy}$
implies the event:
\begin{eqnarray*}
B_{xy}(i,j) &:=& \{\mbox{no visits to }y\mbox{ before }\tau_x(i)\\
&&\hspace{30mm}
\mbox{ and not less than }j\mbox{ visits to }y\mbox{ after }\tau_x(i)\}.
\end{eqnarray*}
Altogether these imply the following upper bound
\begin{eqnarray}\label{estim.2}
\Prob\{B(x,y,i,j)\} &\le& 
\Prob\{B\cap B_{xyx}\}+\Prob\{B\cap B_{yxy}\}\nonumber\\
&&\hspace{10mm}+\Prob\{B\cap B_{xy}(i,j)\}+\Prob\{B\cap B_{yx}(j,i)\}.
\end{eqnarray}
Let us estimate every probability on the right hand side here. 
Since $\tau_x(i)$ is a Markov time, 
\begin{eqnarray*}
\lefteqn{\Prob\{B\cap B_{xy}(i,j)\}}\\ 
&=& \sum_{k\le n}\Prob\{\tau_x(i)=k,
\mbox{not less than }j\mbox{ visits to }y
\mbox{ within time interval }[k+1,n]\}\\
&=& \sum_{k\le n}\Prob\{\tau_x(i)=k\} \Prob\{l(n-k,y-x)\ge j\}\\
&\le& \sum_{k\le n}\Prob\{\tau_x(i)=k\} \Prob\{l(n,y-x)\ge j\},
\end{eqnarray*}
because the event $\{l(n,y-x)\ge j\}$ can only increase as $n$ grows.
Therefore,
\begin{eqnarray*}
\Prob\{B\cap B_{xy}(i,j)\} &\le& 
\Prob\{\tau_x(i)\le n\} \Prob\{l(n,y-x)\ge j\}\\
&=& \Prob\{l(n,x)\ge i\} \Prob\{l(n,y-x)\ge j\}.
\end{eqnarray*}
Then summation over all $x\not= y$ implies that
\begin{eqnarray*}
\lefteqn{\sum_{x\not=y}\bigl(\Prob\{B\cap B_{xy}(i,j)\}
+\Prob\{B\cap B_{yx}(j,i)\}\bigr)}\\
&\le& \sum_{x\not=y}\Bigl(\Prob\{l(n,x)\ge i\} \Prob\{l(n,y-x)\ge j\}
+\Prob\{l(n,y)\ge j\} \Prob\{l(n,x-y)\ge i\}\Bigr)\\
&\le& \sum_{x,y}\Prob\{l(n,x)\ge i\} \Prob\{l(n,y)\ge j\}.
\end{eqnarray*}
Together with non-negativity of increments of the function $f$ 
it implies that
\begin{eqnarray}\label{Bxy}
\lefteqn{\sum_{i,j=1}^n \Delta f(i)\Delta f(j) \biggl(
\sum_{x\not=y}\bigl(\Prob\{B\cap B_{xy}(i,j)\}+\Prob\{B\cap B_{yx}(j,i)\}\bigr)}
\nonumber\\
&&\hspace{25mm}-\sum_{x,y}\Prob\{l(n,x)\ge i\}\Prob\{l(n,y)\ge j\}\biggr)
\ \le\ 0.
\end{eqnarray}
Further, the event $B_{xyx}$ may be described as follows:
firstly the site $x$ is visited at least once, say $t\ge 1$ times, 
then the site $y$ is visited one or more times, say $s\ge 1$ times, 
and then again the site $x$ is visited, which is followed by 
visits to $x$ and $y$ in an arbitrary order. Thus, for $i\ge j$,
\begin{eqnarray*}
\lefteqn{\Prob\{B\cap B_{xyx}\}}\\
&\le& \sum_{t,s,k_1<\ldots<k_{t+s+1}\le n}
\Prob\{S_{k_1}=\ldots=S_{k_t}=x,
S_{k_{t+1}}=\ldots=S_{k_{t+s}}=y,\\
&&\hspace{15mm} S_{k_{t+s+1}}=x,
\mbox{ there are no other visits to }x\mbox{ and }y\mbox{ up to }k_{t+s+1},\\
&&\hspace{45mm}
S_k=x\mbox{ at least }i-t-1\mbox{ times past }k_{t+s+1}\}.
\end{eqnarray*}
and similarly for $j\ge i$
\begin{eqnarray*}
\lefteqn{\Prob\{B\cap B_{xyx}\}}\\
 &\le& \sum_{t,s,k_1<\ldots<k_{t+s+1}\le n}
\Prob\{S_{k_1}=\ldots=S_{k_t}=x,
S_{k_{t+1}}=\ldots=S_{k_{t+s}}=y,\\
&&\hspace{10mm} S_{k_{t+s+1}}=x,
\mbox{ there are no other visits to }x\mbox{ and }y\mbox{ up to }k_{t+s+1},\\
&&\hspace{45mm}
S_k=y\mbox{ at least }j-s\mbox{ times past }k_{t+s+1}\}.
\end{eqnarray*}
Summing up for all $x$ and $y$ 
we arrive at the following upper bound
\begin{eqnarray*}
\lefteqn{\sum_{x\not=y}\Prob\{B\cap B_{xyx}\}}\\
&\le& \Prob^{i-1}\{\tau<\infty\}
\sum_{z\in\Z^d,\ r_1<r_2<r_3\le n}
\Prob\{S_{r_2}-S_{r_1}=z,S_{r_3}-S_{r_2}=-z\}\\
&=& (1-\gamma)^{i-1} \sum_{z\in\Z^d,\ r_1<r_2<r_3\le n}
\Prob\{S_{r_2}-S_{r_1}=z,S_{r_3}-S_{r_2}=-z\}
\end{eqnarray*}
in the case $i\ge j$ and similarly with coefficient 
$(1-\gamma)^{j-1}$ in the case $j\ge i$. Since
$$
\sum_{z\in\Z^d} \Prob\{S_{r_2}-S_{r_1}=z,S_{r_3}-S_{r_2}=-z\}
\ =\ \Prob\{S_{r_3}-S_{r_1}=0\},
$$
we get, for $i\ge j$,
\begin{eqnarray*}
\sum_{x\not=y}\Prob\{B\cap B_{xyx}\} &\le& 
(1-\gamma)^{i-1} \sum_{r_1<r_3\le n} (r_3-r_1)\Prob\{S_{r_3}-S_{r_1}=0\}\\
&=& (1-\gamma)^{i-1} \sum_{r=1}^n r(n-r)\Prob\{S_r=0\},
\end{eqnarray*}
which together with (\ref{upper.1}), (\ref{estim.2}) and (\ref{Bxy})
shows that the variance of $G_n(f)$ does not exceed
\begin{eqnarray*}
\E G_n(f^2)
+4\sum_{j\le i,\ i,j=1}^n \Delta f(i)\Delta f(j)
(1-\gamma)^{i-1} \sum_{r=1}^n r(n-r)\Prob\{S_r=0\}.
\end{eqnarray*}
The sum of $\Delta f(j)$ from $j=1$ to $i$ equals $f(i)$,
hence the desired upper bound for $\Var G_n(f)$.
\qed

\section{Proof of Theorem \ref{t1.2}}
\label{sec:proof1}

Without loss of generality we assume $f(0)=0$.
Any function $f:\Z^+\to\R$ with $f(0)=0$ is decomposable into a
difference of two non-decreasing functions, $f=f_1-f_2$, where
\begin{eqnarray}\label{f=f1-f2}
f_1(j)\ =\ \sum_{i=1}^j (f(i)-f(i-1))^+,\quad
f_2(j)\ =\ \sum_{i=1}^j (f(i)-f(i-1))^-.
\end{eqnarray}
Since
\begin{eqnarray*}
f_k(j) &\le& \sum_{i=1}^j |f(i)-f(i-1)|
\ \le\ 2\sum_{i=1}^j |f(i)|,\quad k=1,\ 2,
\end{eqnarray*}
we get the following upper bound
\begin{eqnarray*}
\sum_{j=1}^\infty f_k^2(j)(1-\gamma)^j &\le& 
4\sum_{j=1}^\infty (1-\gamma)^j\biggl(\sum_{i=1}^j |f(i)|\biggr)^2\\
&\le& 4\sum_{j=1}^\infty (1-\gamma)^j j\sum_{i=1}^j f^2(i)\\
&=& 4\sum_{i=1}^\infty  f^2(i) \sum_{j=i}^\infty (1-\gamma)^j j
\ \le\ \frac{4}{\gamma^2}\sum_{i=1}^\infty  f^2(i)(i+1)(1-\gamma)^i.
\end{eqnarray*}
Therefore, the condition (\ref{cond.on.f}) implies that
\begin{eqnarray}\label{cond.on.f.deco}
\sum_{j=1}^\infty f_k^2(j)(1-\gamma)^j &<& \infty,\quad k=1,\ 2.
\end{eqnarray}
Hence, without loss of generality we assume that $f$ is
a non-decreasing function satisfying (\ref{cond.on.f.deco}) and $f(0)=0$.

The transience of the random walk $(S_n)_{n\ge 0}$ is equivalent to
the convergence of the series
\begin{eqnarray}\label{f3.2}
\sum_{n=1}^\infty\Prob\{S_n=0\} &<& \infty.
\end{eqnarray}
The condition (\ref{cond.on.f.deco}) allows us to apply the upper bound
(\ref{upper.for.Q}) to $f^2$ and to conclude that 
$\E G_n(f^2)\le c_1n$ for some $c_1<\infty$.
Since $f$ is non-decreasing and $f(0)=0$, 
$\Delta f(i)\le f(i)$. Therefore, by Lemma \ref{l:var}, 
\begin{eqnarray*}
\Var G_n(f) &\le& c_1n+c_2n\sum_{i=1}^n f^2(i)
(1-\gamma)^{i-1} \sum_{r=1}^n r\Prob\{S_r=0\}\\
&\le& c_1n+c_3n \sum_{r=1}^n r\Prob\{S_r=0\},
\end{eqnarray*}
again by the condition (\ref{cond.on.f.deco}). In view of (\ref{f3.2}),
$$
a_n\ :=\ \frac{1}{n}\sum_{r=1}^n r\Prob\{S_r=0\}\ \to 0
\quad\mbox{as }n\to\infty
$$
and hence 
\begin{eqnarray}\label{f3.5}
\Var\frac{G_n(f)}{n} &\le& \frac{c_1}{n}+c_3a_n\to 0\quad \mbox{as }n\to\infty,
\end{eqnarray}
which is equivalent to the convergence $(G_n(f)-\E G_n(f))/n\to 0$ in $L_2$.
Together with the convergence (\ref{Gf.conv}) this completes
the proof of $L_2$-convergence stated in Theorem \ref{t1.2}.

For the proof of the almost sure convergence, first let us notice
that (\ref{f3.2}) yields
\begin{eqnarray*}
\sum_{n=1}^\infty {\frac{a_n}{n}} &=& \sum_{n=1}^\infty
{\frac{1}{n^2}{\sum_{r=1}^n r\Prob\{S_r=0\}}}\\
&\le& \sum_{r=1}^\infty \Prob\{S_r=0\} 
\ r\sum_{n=r}^\infty \frac{1}{n^2}\ \le\ 
2\sum_{r=1}^\infty \Prob\{S_r=0\}\ <\ \infty.
\end{eqnarray*}
Hence we can apply Lemma \ref{l3.1} proven below to 
the sequence $\{a_n\}_{n\ge 1}$, so, for any fixed $\delta>0$, 
there is an increasing subsequence $\{n_r\}_{r\ge 1}$ such that
$\sum_{r=1}^\infty a_{n_r}<\infty$ and
$\sqrt{1+\delta}n_{r-1}\le n_{r+1}\le(1+\delta)n_r$ for all $r$.
 
Using Chebyshev's inequality, the upper bound (\ref{f3.5}) and
the convergence (\ref{Gf.conv}) we conclude that, for any $\varepsilon>0$,
\begin{eqnarray*}
\sum_{r=1}^\infty
\Prob\left\{\Bigl|\frac{G_{n_r}(f)}{\E G_{n_r}(f)}-1\Bigl|
>\varepsilon\right\} &\le& \sum_{r=1}^\infty
\frac{\Var G_{n_r}(f)}{{\varepsilon}^2 {(\E G_{n_r}}(f))^2}\\
&\le& \frac{C}{\varepsilon^2}\sum_{r=1}^\infty (1/n_r+a_{n_r})<\infty.
\end{eqnarray*}
Then it follows from the Borel--Cantelli lemma that
\begin{equation}\label{f3.16}
\frac{G_{n_r}(f)}{\E G_{n_r}(f)}\stackrel{a.s.}\to 1
\quad \mbox{as }r\to\infty.
\end{equation}
Further, for any $n$ there exists $r$ such that $n_r\le n\le n_{r+1}$ 
and, hence
\begin{equation}\label{f3.17}
\frac{G_{n_r}(f)}{\E G_{n_{r+1}}(f)}\ \le\
\frac{G_n(f)}{\E G_n(f)}\ \le\
\frac{G_{n_{r+1}}(f)}{\E G_{n_r}(f)}.
\end{equation}
It follows from (\ref{Gf.conv}) that
$$
\frac{\E G_{n_{r+1}}(f)}{\E G_{n_r}(f)}\ \sim\
\frac{n_{r+1}}{n_r}\quad\mbox{as }r\to\infty.
$$
Moreover, $n_r<n_{r+1}\le(1+\delta)n_r$ for all $r$. 
Then (\ref{f3.16}) and (\ref{f3.17}) imply that
\begin{eqnarray*}
\frac{1}{1+\delta}\ \le\
\liminf_{n\to\infty}{\frac{G_n(f)}{\E G_n(f)}} &\le&
\limsup_{n\to\infty}{\frac{G_n(f)}{\E G_n(f)}}\ \le\ 1+\delta
\quad \mbox{a.s.}
\end{eqnarray*}
Due to arbitrary choice of $\delta>0$,
the a.s. convergence $G_n(f)/\E G_n(f)\to 1$ follows.
\qed

In the last proof, we have made use of the following auxiliary result.

\begin{Lemma}\label{l3.1}
Let $v_n\ge 0$ and $\sum_{n=1}^\infty \frac{v_n}{n}<\infty$. 
Then, for any fixed $\delta>0$, there exists an increasing
subsequence $\{n_r\}_{r\ge 1}$ such that
$\sum_{r=1}^\infty v_{n_r}<\infty$ 
and $\sqrt{1+\delta}n_{r-1}\le n_{r+1}\le (1+\delta)n_r$ for all $r\ge 1$.
\end{Lemma}

\proof
Let us fix an arbitrary $b\in(1,2)$ and identify a $K=K(b)$ such that 
$[b^K]-[b^{K-1}]\ge 2$. For $r\ge 1$, choose
$$
n_r\in[[b^{K+r-2}]+1,[b^{K+r-1}]]\mbox{ such that }
v_{n_r}\ =\ \min_{[b^{K+r-2}]+1\le n \le[b^{K+r-1}]}{v_n}.
$$
By this construction, 
\begin{eqnarray*}
\sum_{n=1}^\infty \frac{v_n}{n} &\ge&
v_{n_1}\sum_{n=1}^{[b^K]} \frac{1}{n}+
v_{n_2}\sum_{n=[b^K]+1}^{[b^{K+1}]}\frac{1}{n}+\ldots+
v_{n_r}\sum_{n=[b^{K+r-2}]+1}^{[b^{K+r-1}]}\frac{1}{n}+\ldots
\end{eqnarray*}
Since
\begin{eqnarray*}
\sum_{n=\left[ b^{K+r-2} \right]+1}^{\left[ b^{K+r-1}\right]}\frac{1}{n}
&\to& \log b\ >\ 0\quad\mbox{as } r\to\infty,
\end{eqnarray*}
the convergence of the series $\sum_n\frac{v_n}{n}$ 
guarantees convergence of the series $\sum_r v_{n_r}$. 
Also, for all $r$,
\begin{eqnarray*}
\frac{n_{r+1}}{n_{r-1}}\ \ge\ \frac{b^{K+r-1}}{b^{K+r-2}}\ =\ b
\quad\mbox{ and }\quad
\frac{n_{r+1}}{n_r}\ \le\ \frac{b^{K+r}}{b^{K+r-2}}\ =\ b^2,
\end{eqnarray*}
so the lemma conclusion follows if we take $b=\sqrt{1+\delta}$.
\qed

\section{Proof of Theorem \ref{thm:exp.conv}}
\label{sec:max}

To prove Theorem \ref{thm:exp.conv}, let us first consider 
the maximal local time, $l(n):=\max\{l(n,x),\ x\in\Z^d\}$.
Theorem 13 in Erd\H{o}s and Taylor \cite{ET}
states a strong limit theorem for $l(n)$: 
for a simple random walk in $d\ge 3$ dimensions,
\begin{eqnarray}\label{thm:max.l.t}
\frac{l(n)}{\log n} &\to&
\frac{1}{\log\frac{1}{1-\gamma}}\ =:\ \frac{1}{\lambda_*}\quad\mbox{as }n\to\infty
\mbox{ with probability }1.
\end{eqnarray}
The proof in \cite{ET} is split into two parts, 
dealing with upper and lower bounds. There is some issue with 
the proof of the upper bound, that is,
\begin{eqnarray}\label{thm:max.l.t.upper}
\limsup_{n\to\infty}\frac{l(n)}{\log n} &\le& \frac{1}{\lambda_*}
\quad\mbox{with probability }1.
\end{eqnarray}
The proof suggested in \cite{ET} is based on the inequality for tails
\begin{eqnarray}\label{l.n.tail.geom}
\Prob\{l(n)>t\} &\le& n\Prob\{l(n,0)>t\}\quad\mbox{for all }t>0,
\end{eqnarray}
and on the observation that the number of returns to the origin,
$l(n,0)$, is dominated by a geometrically distributed random variable 
with parameter $1-\gamma$. Notice that justification of (\ref{l.n.tail.geom}) 
in \cite{ET} is not complete because it is based there on the assumption 
that all sites visited by the random walk---clearly 
not more than $n$---can be treated in the same way as the origin.
This point requires further justification 
because the set of visited sites is random.
The same issue occurs is the proof of Theorem 1 in Revesz \cite{Revesz}.
Notice that this set is contained in the ball 
of radius $n$, which leads to the coefficient $n^d$ instead of $n$
on the right hand side of (\ref{l.n.tail.geom}) which in its turn
leads to the constant $d/\lambda_*$ on the right hand side of 
(\ref{thm:max.l.t.upper}) instead of $1/\lambda_*$.

The last issue may be resolved in different ways, particularly,
we may condition on the non-zero value of $S_1$,
\begin{eqnarray*}
\Prob\{l(n)>t\} &\le& 
\Prob\{l(n,0)>t\}+\sum_{x\not= 0}\Prob\{l(n-1,x)>t\mid S_1=x\}\Prob\{S_1=x\},
\end{eqnarray*}
followed by an induction argument on $n$. 
Hence, the upper bound (\ref{l.n.tail.geom})
holds for any transient random walk in any dimensions. Therefore, 
\begin{eqnarray*}
\Prob\{l(n)>t\} &\le& n\Prob\{l(\infty,0)>t\}\ \le\ n(1-\gamma)^{t-1}
\ =\ ne^{-(t-1)\lambda_*},
\end{eqnarray*}
which implies, for all $\varepsilon>0$ and $m\in\{0,1,2,\ldots\}$,
the following upper bound
\begin{eqnarray*}
\Prob\{C(n,\varepsilon,m)\} &\le& 
\frac{2}{\log n\ldots\log_{(m-2)}n\log_{(m-1)}^{1+\varepsilon}n},
\end{eqnarray*}
for the events
\begin{eqnarray*}
C(n,\varepsilon,m) &:=& \Bigl\{l(2n)>1+
\frac{1}{\lambda_*}(\log n+\ldots+\log_{(m-1)}n
+(1+\varepsilon)\log_{(m)}n)\Bigr\};
\end{eqnarray*}
hereinafter $\log_{(m)}x$ denotes the $m$-fold iterated logarithms, that is, 
$$
\log_{(0)}(x)=x,\quad
\log_{(m)}(x)=\log\log_{(m-1)}x\quad\mbox{for all }m\ge 1.
$$
Therefore, the series $\sum_{k=1}^\infty\Prob\{C(2^k,\varepsilon,m)\}$ converges, 
hence by the Borel-Cantelli lemma, only finitely many
of $C(2^k,\varepsilon,m)$ occur, with probability $1$. 
For any $n\in[2^k,2^{k+1})$ and the event
\begin{eqnarray*}
B(n,\varepsilon,m) &:=& \Bigl\{l(n)>1+
\frac{1}{\lambda_*}(\log n+\ldots+\log_{(m-1)}n
+(1+\varepsilon)\log_{(m)}n)\Bigr\},
\end{eqnarray*}
we have inclusion $B(n,\varepsilon,m)\subseteq C(2^k,\varepsilon,m)$,
and thus only finitely many of $B(2^k,\varepsilon,m)$ occur, with probability $1$.
In other words, we arrive at the following result.

\begin{Proposition}\label{l:ln.limit.thm}
For all $\varepsilon>0$ and $m\in\{0,1,2,\ldots\}$, 
\begin{eqnarray*}
l(n) &\le& 
\frac{1}{\lambda_*}\Bigl(\log n+\ldots+\log_{(m-1)}n
+(1+\varepsilon)\log_{(m)}n\Bigl),
\end{eqnarray*}
for all $n\ge N$ where $N$ is finite with probability $1$.
\end{Proposition}

Notice that, for a simple random walk in $\Z^d$, $d\ge 3$, an upper a.s.\ bound
$\lambda_*^{-1}(\log n+(1+\varepsilon)\log\log n)$ and---in the case 
of $d\ge 4$---a lower a.s.\ bound 
$\lambda_*^{-1}(\log n-(3+\varepsilon)\log\log n)$ 
is derived by Revesz in \cite{Revesz} following a different technique;
he has also proved that 
$l(n)\ge\lambda_*^{-1}(\log n+(1-2/(d-2)-\varepsilon)\log\log n)$ 
infinitely often a.s.
The maximal local time for a zero drift random walk on $\Z$ 
with finite variance---which is clearly recurrent---was studied
by Kesten in \cite{Kesten}.

Let us proceed with the proof of Theorem \ref{thm:exp.conv},
we start with the case (ii).
Introducing two non-decreasing functions, $f_1$ and $f_2$
as in (\ref{f=f1-f2}), we notice that then 
$f_k(i)\le\widetilde C e^{i\lambda_*}/i\log^{2+\varepsilon}i$ 
for $k=1$, $2$, because
\begin{eqnarray*}
\sum_{i=1}^n \frac{e^{i\lambda_*}}{i\log^{2+\varepsilon}i} 
&=& O\Bigl(\frac{e^{n\lambda_*}}{n\log^{2+\varepsilon}n}\Bigr)
\quad\mbox{as }n\to\infty.
\end{eqnarray*}
Hence, without loss of generality we assume that $f$ 
is non-decreasing with $f(0)=0$.

The function $f$ satisfies the condition (\ref{cond.on.f.w}), but not 
(\ref{cond.on.f.deco}), and this generates a certain difficulty 
we need to overcome.
To this end, let us introduce two sequences of truncated functions
\begin{eqnarray*}
f_n(i) &:=& f(i)\I\bigl\{i\le \lambda_*^{-1}\log n\bigr\},\\
f_n^+(i) &:=& f(i)\I\bigl\{\lambda_*^{-1} \log n<i\le 
\lambda_*^{-1}(\log n+b_n)\bigr\},
\end{eqnarray*}
where $b_n=\log_{(2)}n+2\log_{(3)}n$ 
and make use of the following decomposition
\begin{eqnarray*}
G_n(f) &=& G_n(f_n)+G_n(f_n^+)+G_n(f-f_n-f_n^+).
\end{eqnarray*}
Since $f$ satisfies the condition (\ref{cond.on.f.w}), we get equivalences
\begin{eqnarray}\label{Gn.equiv.n}
\E G_n(f_n) &\sim& \E G_n(f)\ \sim\ n\gamma^2
\sum_{j=1}^\infty f(j)(1-\gamma)^{j-1}\quad\mbox{as }n\to\infty,
\end{eqnarray}
and, by (\ref{upper.for.Q}), the following upper bound
\begin{eqnarray*}
\E G_n(f_n^2) &\le& c_1\frac{n}{\log n\log_{(2)}^{2+\varepsilon}n}
\E G_n(f_n)\ \le\ c_2\frac{n^2}{\log n\log_{(2)}^2 n}.
\end{eqnarray*}
Further, it follows from the condition (\ref{bound.for.P.Sn.conv}) that
\begin{eqnarray*}
\sum_{r=1}^n r(n-r)\Prob\{S_r=0\} &\le& n\sum_{r=1}^n r\Prob\{S_r=0\}\ \le\ c_3n\log n.
\end{eqnarray*}
In addition,
\begin{eqnarray*}
\sum_{i=1}^n f_n(i)\Delta f_n(i)(1-\gamma)^{i-1} &\le& 
\sum_{i=1}^{\lambda_*^{-1}\log n} f^2(i)e^{-(i-1)\lambda_*}\\
&\le& c_4\frac{e^{i\lambda_*}}{i^2\log^4 i}\Big|_{i=\lambda_*^{-1}\log n}\\
&\le& c_5\frac{n}{\log^2 n\log_{(2)}^4 n}.
\end{eqnarray*}
Substituting the last three bounds into the right hand side
of the inequality provided by Lemma \ref{l:var}, we derive that
\begin{eqnarray*}
\Var G_n(f_n) &\le& c_6n^2/\log n\log_{(2)}^2n.
\end{eqnarray*}
Choose a subsequence $n_r=[e^{r/\log^{1/4}r}]$, then, by Chebyshev's inequality, 
the last upper bound and (\ref{Gn.equiv.n}), we conclude that
\begin{eqnarray*}
\sum_{r=1}^\infty
\Prob\left\{\Bigl|\frac{G_{n_r}(f_{n_r})}{\E G_{n_r}(f_{n_r})}-1\Bigl|
>\frac{1}{\log^{1/8} r}\right\} &\le& \sum_{r=1}^\infty
\frac{\log^{1/4} r\Var G_{n_r}(f_{n_r})}{{(\E G_{n_r}}(f_{n_r}))^2}\\
&\le& c_7\sum_{r=1}^\infty \frac{1}{r\log^{3/2}r}\ <\ \infty,
\end{eqnarray*}
which allows us to apply the Borel--Cantelli lemma, hence obtaining
\begin{eqnarray*}
\frac{G_{n_r}(f_{n_r})}{\E G_{n_r}(f_{n_r})} 
&\stackrel{a.s.}\to& 1 \quad \mbox{as }r\to\infty.
\end{eqnarray*}
Similar to (\ref{f3.17}), if $n_r\le n\le n_{r+1}$ then
\begin{eqnarray*}
\frac{G_{n_r}(f_{n_r})}{\E G_{n_{r+1}}(f_{n_{r+1}})} &\le&
\frac{G_n(f_n)}{\E G_n(f_n)}\ \le\
\frac{G_{n_{r+1}}(f_{n_{r+1}})}{\E G_{n_r}(f_{n_r})}.
\end{eqnarray*}
In addition, by (\ref{Gn.equiv.n}),
$$
\frac{\E G_{n_{r+1}}(f_{n_{r+1}})}{\E G_{n_r}(f_{n_r})}\ \sim\
\frac{n_{r+1}}{n_r}\ \sim\
e^{\frac{r+1}{\log^{1/4}(r+1)}-\frac{r}{\log^{1/4}r}}\ \to\ 1\quad\mbox{as }r\to\infty.
$$
Therefore,
\begin{equation}\label{f3.16.e}
\frac{G_n(f_n)}{\E G_n(f_n)}\stackrel{a.s.}\to 1
\quad \mbox{as }n\to\infty.
\end{equation}
Further, it follows from (\ref{upper.for.Q}) that, for all $m\le n$,
\begin{eqnarray*}
\E\frac{G_n(f_m^+)}{n} &\le& 
\sum_{j=[\lambda_*^{-1}\log m]}^{[\lambda_*^{-1}(\log m+b_m)]}
f(j)(1-\gamma)^{j-1}\\
&\le& c_7
\sum_{j=[\lambda_*^{-1}\log m]}^{[\lambda_*^{-1}(\log m+b_m)]}
\frac{1}{{j\log^{2+\varepsilon}j}}\\
&=& O\biggl(\frac{b_m}{\log m\log_{(2)}^{2+\varepsilon}m}\biggr)
\ =\ O\biggl(\frac{1}{\log n\log_{(2)}^{1+\varepsilon}n}\biggr),
\end{eqnarray*}
if $m\ge n/2$.
Applying Chebyshev's inequality to non-negative random variables
$G_{n_{k+1}}(f_{n_k}^+)/n_{k+1}$ with $n_k=2^k$, 
we get the following series convergence
\begin{eqnarray*}
\sum_{k=1}^\infty
\Prob\Bigl\{\frac{G_{n_{k+1}}(f_{n_k}^+)}{n_{k+1}}\ge
\frac{1}{\log^{\varepsilon/2}k}\Bigr\} 
&\le& \sum_{k=1}^\infty \frac{c_8/k\log^{1+\varepsilon}k}
{1/\log^{\varepsilon/2}k}
\ <\ \infty,
\end{eqnarray*}
which in its turn implies by the Borel--Cantelli lemma that
\begin{eqnarray*}
\frac{G_{n_{k+1}}(f_{n_k}^+)}{n_{k+1}} &\to& 0
\quad\mbox{as }k\to\infty\mbox{ with probability }1.
\end{eqnarray*}
In addition, for $n_k\le n\le n_{k+1}$,
\begin{eqnarray*}
\frac{G_n(f_n^+)}{n} &\le& \frac{G_{n_{k+1}}(f_{n_k}^+)}{n_k}
\ =\ 2\frac{G_{n_{k+1}}(f_{n_k}^+)}{n_{k+1}},
\end{eqnarray*}
hence
\begin{eqnarray}\label{f3.16.e.+}
\frac{G_n(f_n^+)}{n} &\to& 0
\quad\mbox{as }n\to\infty\mbox{ with probability }1.
\end{eqnarray}
Finally, since $l(n)$ is the largest local time,
\begin{eqnarray*}
\{G_n(f-f_n-f_n^+)>0\} &\subseteq& 
\{l(n)>\lambda_*^{-1}\bigl(\log n+\log_{(2)}n+2\log_{(3)}n\bigr)\},
\end{eqnarray*}
which implies a.e. convergence $G_n(f-f_n-f_n^+)\to 0$ as $n\to\infty$,
due to Proposition \ref{l:ln.limit.thm} with $m=3$.
Together with (\ref{f3.16.e}), (\ref{f3.16.e.+}) and (\ref{Gn.equiv.n}) 
it implies the desired convergence (\ref{Gf.conv.L21}) in the case (ii).

In the case (i), the proof requires some alterations.
Consider a sequence of truncated functions
\begin{eqnarray*}
f_n(i) &:=& f(i)\I\biggl\{i\le \frac{\lambda_*^{-1}\eta}{2}\log n\biggr\},
\end{eqnarray*}
and make use of the following decomposition
\begin{eqnarray*}
G_n(f) &=& G_n(f_n)+G_n(f-f_n).
\end{eqnarray*}
As above, the equivalences (\ref{Gn.equiv.n}) hold and, by (\ref{upper.for.Q}),
\begin{eqnarray*}
\E G_n(f_n^2) &\le& c_9n^{\eta/2}\E G_n(f_n)\ \le\ c_{10}n^{1+\eta/2}.
\end{eqnarray*}
Further, it follows from the condition (\ref{bound.for.P.Sn=0}) that
\begin{eqnarray*}
\sum_{r=1}^n r(n-r)\Prob\{S_r=0\} &\le& n\sum_{r=1}^n r\Prob\{S_r=0\}
\ \le\ c_{11}n^{2-\eta}.
\end{eqnarray*}
In addition,
\begin{eqnarray*}
\sum_{i=1}^n f_n(i)\Delta f_n(i)(1-\gamma)^{i-1} &\le& 
\sum_{i=1}^{\frac{\lambda_*^{-1}\eta}{2}\log n} f^2(i)e^{-(i-1)\lambda_*}
\ \le\ c_{12}n^{\eta/2}.
\end{eqnarray*}
Substituting the last three bounds into the right hand side
of the inequality provided by Lemma \ref{l:var}, we derive that
\begin{eqnarray*}
\Var G_n(f_n) &\le& c_{13}n^{2-\eta/2},
\end{eqnarray*}
since $\eta<1$.
Then similar to the case (ii) we deduce (\ref{f3.16.e}).
Further, it follows from (\ref{upper.for.Q}) that, for all $m\le n$,
\begin{eqnarray*}
\E\frac{G_n(f-f_m)}{n} &\le& 
\sum_{j=[\frac{\lambda_*^{-1}\eta}{2}\log m]}^n f(j)(1-\gamma)^{j-1}\\
&\le& c_{14} \sum_{j=[\frac{\lambda_*^{-1}\eta}{2}\log m]}^\infty
\frac{1}{j^{2+\varepsilon}}
\ =\ O\biggl(\frac{1}{\log^{1+\varepsilon}n}\biggr),
\end{eqnarray*}
if $m\ge n/2$. Again similar to the case (ii), we deduce from the last bound that
\begin{eqnarray*}
\frac{G_n(f-f_n)}{n} &\to& 0
\quad\mbox{as }n\to\infty\mbox{ with probability }1,
\end{eqnarray*}
which together with (\ref{f3.16.e}) and equality 
$G_n(f)=G_n(f_n)+G_n(f-f_n)$ implies (\ref{Gf.conv.L21}) in the case (i). 
The proof of Theorem \ref{thm:exp.conv} is complete.
\qed

\section*{Acknowledgments}

The authors are very thankful to the referee for valuable comments.

\end{document}